\documentclass[11pt]{article}
\usepackage{latexsym,amssymb}
\def\demo{\noindent{\bf Proof. }}
\def\QED{\hfill$\Box$}
\newtheorem{Theorem}{Theorem}[section]

\newtheorem{Corollary}[Theorem]{Corollary}

\begin{document}

\topmargin3mm
\hoffset=-1cm
\voffset=-1.5cm
\begin{center}

{\Large\bf Unmixed bipartite graphs}\\

\vspace{6mm}


\addtocounter{footnote}{1}


\medskip

Rafael
H. Villarreal\footnote{Partially supported
by CONACyT grant 49251-F and SNI, M\'exico.}\\

{\small Departamento de
Matem\'aticas}\vspace{-1mm}\\
{\small Centro de Investigaci\'on y de Estudios
Avanzados del
IPN}\vspace{-1mm}\\
{\small Apartado Postal
14--740}\vspace{-1mm}\\
{\small 07000 Mexico City, D.F.}\vspace{-1mm}\\
{\small
e-mail: {\tt
vila@math.cinvestav.mx}}\vspace{4mm}

\end{center}

\date{}

\begin{abstract}
\noindent In this note we give a combinatorial characterization of
all the  
unmixed bipartite graphs. 
\end{abstract}

\centerline{\small \bf Resumen}\vspace{-2mm} 
\begin{verse}\small En esta nota nosotros presentamos una 
caracterizaci\'on combinatoria de \\ todas las gr\'aficas bipartitas 
no-mezcladas.
\end{verse}

\section{Unmixed graphs}

In the sequel we use \cite{Har} as a reference for 
standard terminology and notation on graph theory. 

\medskip

Let $G$ be a simple graph with vertex set $V(G)$ and edge set $E(G)$. 
A subset $C\subset V(G)$ 
is a {\it minimal vertex cover\/} of 
$G$ if: (1) every edge of $G$ is incident with one vertex in $C$, 
and (2) there is no proper subset of $C$ with the first 
property. If $C$ satisfies condition (1) only, then $C$ is 
called a {\it vertex cover\/} of $G$. Notice that $C$ is a 
minimal vertex cover if and only if $V(G)\setminus C$ is a 
maximal independent set. 
A graph $G$ is called {\it unmixed\/} if all the minimal vertex
covers   
of $G$ have the same number of elements and it is called 
{\it well covered\/} \cite{plummer-unmixed} if all the maximal
independent sets of $G$ have the same number of elements. 

\medskip

The notion of unmixed graph is related to some other graph 
theoretical and algebraic properties. 
The following implications hold for any graph without isolated 
vertices \cite{berge1,Har,Vi2}:
$$
\begin{array}{ccccccc}
\mbox{Cohen-Macaulay}&\Longrightarrow& \mbox{unmixed}&
\Longrightarrow&B\mbox{-graph} &\Longrightarrow
&\mbox{vertex-critical.} \\ 
\end{array}
$$ 
Structural aspects of Cohen-Macaulay bipartite graphs were first
studied in \cite{EV}. In loc. cit. it is shown that $G$ is
Cohen-Macaulay if and only if the simplicial complex $\Delta_G$
generated by the maximal independent sets of $G$ is shellable. 
The main result that we present in
this note is the following combinatorial characterization of all the  
unmixed bipartite graphs. Our result is inspired by a criterion of 
 Herzog and Hibi \cite[Theorem~3.4]{herzog-hibi} 
that describe all Cohen-Macaulay bipartite graphs
in combinatorial terms. 

\begin{Theorem}\label{unmixed-bip} Let $G$ be a bipartite graph
without isolated vertices. Then $G$ is unmixed if and only 
if there is a bipartition $V_1=\{x_1,\ldots,x_g\}$,
$V_2=\{y_1,\ldots,y_g\}$ of $G$ such that: 
{\rm (a)} $\{x_i,y_i\}\in E(G)$ for all $i$, and {\rm (b)} 
if $\{x_i,y_j\}$ and $\{x_j,y_k\}$ are in $E(G)$ and $i,j,k$ are
distinct, then $\{x_i,y_k\}\in E(G)$.
\end{Theorem} 
\demo $\Rightarrow$) Since $G$ is bipartite, 
there is a bipartition $(V_1,V_2)$ of $G$, i.e., 
$V(G)=V_1\cup V_2$, $V_1\cap V_2=\emptyset$, and every edge of $G$ 
joins $V_1$ with $V_2$. Let $g$ be the vertex covering number of $G$, 
i.e., $g$ is the number of elements in any minimal vertex cover 
of $G$. Notice that $V_1$ and $V_2$ are both minimal vertex covers 
of $G$, hence $g=|V_1|=|V_2|$. By K\"onig theorem \cite[Theorem~10.2,
p.~96]{Har} $g$ is 
the maximum number of independent edges of $G$. Therefore after
permutation  
of the vertices we obtain that $V_1=\{x_1,\ldots,x_g\}$, 
$V_2=\{y_1,\ldots,y_g\}$, and that $\{x_i,y_i\}\in E(G)$ for 
$i=1,\ldots,g$. Thus we have proved that (a) holds. To prove (b) 
take $\{x_i,y_j\}$ and $\{x_j,y_k\}$ in $E(G)$ such that $i,j,k$ are
distinct. Assume that $x_i$ is not adjacent to $y_k$. Then there is a
maximal independent set of vertices $A$ containing $x_i$ and $y_k$.
Notice that $|A|=g$ because $G$ is unmixed. Hence $C=V(G)\setminus A$
is a minimal vertex cover of $G$ with $g$ vertices. Since $x_i$ and
$y_k$ are not on $C$, we get that $y_j$ and $x_j$ are both in $C$.
As $C$ intersects $\{x_\ell,y_\ell\}$ in at least one vertex for 
$\ell\neq j$, we obtain that $|C|\geq g+1$, a contradiction. 

$\Leftarrow$) Let $C$ be a minimal vertex cover of $G$. It suffices 
to prove that $C$ intersects $\{x_j,y_j\}$ in exactly one vertex for 
$j=1,\ldots,g$. Assume that $x_j$  and $y_j$ belong to $C$ for some
$j$. If $v\in V(G)$, we denote the {neighbor set}  of $v$ by $N_G(v)$.
Thus there are $x_i\in N_G(y_j)\setminus\{x_j\}$ and $y_k\in
N_G(x_j)\setminus\{y_j\}$ such that $x_i\notin C$ and $y_k\notin C$.
Notice that $i,j,k$ are distinct. Indeed if $i=k$, then $\{x_i,y_i\}$
is an edge of $G$ not covered by $C$, which is impossible. Therefore
using (b) we get that $\{x_i,y_k\}$ is an edge of $C$, a
contradiction. \QED

\medskip

Ravindra \cite{ravindra} has shown a characterization of well covered
bipartite graphs. Namely, $G$ is well covered if
and only if for every edge $\{x,y\}$ in the perfect matching, the
induced subgraph $\langle N_G(x)\cup N_G(y)\rangle$ is a complete
bipartite graph. The advantage of our characterization is that it admits
a natural possible extension to hypergraphs and clutters with a perfect
matching of K\"onig type \cite{morey-reyes-villarreal}.  

\medskip

As a consequence of Theorem~\ref{unmixed-bip} we recover the
following result  
on the structure of unmixed trees. 

\begin{Corollary}{\rm \cite[Theorem~2.4, Corollary~2.5]{Vi2}} 
Let $G$ be a tree with at least three vertices. 
Then $G$ is unmixed if and only 
if there is a bipartition $V_1=\{x_1,\ldots,x_g\}$,
$V_2=\{y_1,\ldots,y_g\}$ of $G$ such that: 
{\rm (a)} $\{x_i,y_i\}\in E(G)$ for all $i$, and {\rm (b)}  
for each $i$ either $\deg(x_i)=1$ or $\deg(y_i)=1$.
\end{Corollary}

\bibliographystyle{plain}

\begin{thebibliography}{99}

\bibitem{berge1}{C. Berge, Some common properties for 
regularizable graphs, edge-critical graphs and $B$-graphs, 
{\it in} ``Theory and practice of combinatorics'' (A. Rosa, 
G. Sabidussi and J. Turgeon, Eds.), North-Holland 
Math. Stud. {\bf 60}, North-Holland, Amsterdam, 1982, 
pp.~31--44. }

\bibitem {EV}{M. Estrada and R. H. Villarreal, Cohen-Macaulay 
bipartite graphs,
Arch. Math. {\bf 68} (1997), 124--128.}

\bibitem{Har}{F. Harary, {\it Graph Theory\/}, Addison-Wesley, 
Reading, MA, 1972.}

\bibitem{herzog-hibi} J. Herzog and T. Hibi, 
Distributive lattices, bipartite graphs and Alexander duality, 
J. Algebraic Combin. {\bf 22} (2005), no. 3, 289--302. 

\bibitem{morey-reyes-villarreal} S. Morey, E. Reyes and R. H.
Villarreal, Cohen-Macaulay, Shellable and unmixed clutters with a
perfect matching of K\"onig type, in progress. 

\bibitem{plummer-unmixed} M. D. Plummer, 
Some covering concepts in graphs, J. Combinatorial Theory {\bf 8}
(1970), 91--98. 

\bibitem{ravindra} G. Ravindra, Well-covered graphs, 
J. Combinatorics Information Syst. Sci. {\bf 2} (1977), no. 1,
20--21.  

\bibitem{Vi2}{R. H. Villarreal, Cohen-{M}acaulay graphs, Manuscripta 
Math. {\bf 66} (1990), 277--293.}

\end{thebibliography}

\end{document}